\documentclass[12pt]{amsart}

\usepackage{amsmath}
\usepackage{amscd}
\usepackage{amssymb}
\usepackage{latexsym}
\usepackage{stmaryrd}
\usepackage{mathbbol}

\input xy
\xyoption{all} \CompileMatrices

\newtheorem{theorem}{Theorem}[section]
\newtheorem*{theorem*}{Theorem}
\newtheorem{lemma}[theorem]{Lemma}
\newtheorem*{lemma*}{Lemma}
\newtheorem{corollary}[theorem]{Corollary}
\newtheorem{proposition}[theorem]{Proposition}

\newtheorem{remark}[theorem]{Remark}
\newtheorem{definition}[theorem]{Definition}

\makeatletter
\def\revddots{\mathinner{\mkern1mu\raise\p@
\vbox{\kern7\p@\hbox{.}}\mkern2mu
\raise4\p@\hbox{.}\mkern2mu\raise7\p@\hbox{.}\mkern1mu}}
\makeatother % ddots gespiegelt
\newcommand{\bgl}{\begin{equation}} %eine Gleichung mit Ziffer
\newcommand{\egl}{\end{equation}}
\newcommand{\bgloz}{\begin{equation*}} %eine Gleichung ohne Ziffer
\newcommand{\egloz}{\end{equation*}}
\newcommand{\bgln}{\begin{eqnarray}} %mehrere Gleichungen mit Ziffer
\newcommand{\egln}{\end{eqnarray}}
\newcommand{\bglnoz}{\begin{eqnarray*}} %mehrere Gleichungen ohne Ziffer
\newcommand{\eglnoz}{\end{eqnarray*}}
\newcommand{\btheo}{\begin{theorem}}
\newcommand{\etheo}{\end{theorem}}
\newcommand{\btheooz}{\begin{theorem*}}
\newcommand{\etheooz}{\end{theorem*}}
\newcommand{\blemma}{\begin{lemma}}
\newcommand{\elemma}{\end{lemma}}
\newcommand{\blemmaoz}{\begin{lemma*}}
\newcommand{\elemmaoz}{\end{lemma*}}
\newcommand{\bproof}{\begin{proof}}
\newcommand{\eproof}{\end{proof}}
\newcommand{\bbew}{\begin{beweis}}
\newcommand{\ebew}{\end{beweis}}
\newcommand{\bremark}{\begin{remark}\em}
\newcommand{\eremark}{\end{remark}}
\newcommand{\bdefin}{\begin{definition}}
\newcommand{\edefin}{\end{definition}}
\newcommand{\bprop}{\begin{proposition}}
\newcommand{\eprop}{\end{proposition}}
\newcommand{\bcor}{\begin{corollary}}
\newcommand{\ecor}{\end{corollary}}
\newcommand{\bfa}{\begin{cases}} %Fallunterscheidung
\newcommand{\efa}{\end{cases}}
\newcommand{\cO}{\mathcal O}

\newcommand{\so}{\text{\tiny{$\cO$}}}

\def\Az{\mathbb{A}}

\def\Qz{\mathbb{Q}}
\def\Rz{\mathbb{R}}

\def\Zz{\mathbb{Z}}

\def\1z{\mathbb{1}}
\newcommand{\fA}{\mathfrak A}

\newcommand{\an}[1]{``#1''} % Anfuehrungsstriche

\newcommand{\ma}{\mapsto} % wird abgebildet auf
\def\SEMI{\mbox{$\times\kern-2pt\vrule height5pt width.6pt \kern3pt $}}

\newcommand{\Aut}{{\rm Aut}\,}

\newcommand{\id}{{\rm id}}

\newcommand{\reg}{^\times} % regulaer
\newcommand{\pos}{_{>0}} % positiv
\newcommand{\dop}{\text{: }} % in Mengen
\newcommand{\falls}{\text{ if }} % bei Fallunterscheidungen
\newcommand{\extalg}{\Lambda^* \,} % aeu{\ss}eres Produkt
\newcommand{\lge}{\left\{} % links geschweift
\newcommand{\rge}{\right\}} % rechts geschweift
\newcommand{\lru}{\left(} % links rund
\newcommand{\rru}{\right)} % rechts rund
\newcommand{\leck}{\left[} % links eckig
\newcommand{\reck}{\right]} % rechts eckig
\newcommand{\lsp}{\left\langle} % links spitz
\newcommand{\rsp}{\right\rangle} % links spitz
\newcommand{\rukl}[1]{\lru #1 \rru} % runde Klammer
\newcommand{\eckl}[1]{\leck #1 \reck} % eckige Klammer
\newcommand{\gekl}[1]{\lge #1 \rge} % geschweifte Klammer
\newcommand{\spkl}[1]{\lsp #1 \rsp} % spitze Klammer
\newcommand{\menge}[2]{\gekl{ #1 \dop #2 }} % Menge
\begin{document}

\title[Erratum]{Erratum to \an{C*-algebras associated with integral domains and crossed products by actions on adele spaces} by J. Cuntz and X. Li}

\author{Joachim Cuntz and Xin Li}

%\begin{abstract}
%We correct the K-theory computation for C*-algebras associated with rings of integers in number fields.
%\end{abstract}

\thanks{Research supported by the Deutsche Forschungsgemeinschaft (SFB 878) and by the ERC through AdG 267079.}

\maketitle

\section{Introduction}
In \cite{Cu-Li}, we had computed the K-theory for C*-algebras
associated with rings of integers in number fields. Unfortunately,
there was a miscalculation in \cite[\S~6.4, case~c)]{Cu-Li}, where
the case of number fields with roots of unity $+1, -1$ and with an
even strictly positive number of real places was treated (i.e. the
case where $\# \gekl{v_{\Rz}} \geq 2$ even). In this case the final
result for the K-theory of the ring C*-algebra $\fA[\so]$ of the
ring of integers $\so$ of our number field should not be
$K_*(\fA[\so]) \cong \extalg (\Gamma) \oplus ((\Zz / 2 \Zz)
\otimes_{\Zz} \extalg (\Gamma))$, but $K_*(\fA[\so]) \cong \extalg
(\Gamma)$. This means that the torsion-free part in \cite[\S~6.4,
case~c)]{Cu-Li} was determined correctly, but the torsion part was
not computed correctly. The correct computation shows that the
K-theory of the ring C*-algebra is torsion-free.

On the whole, the correct final result is the following (compare
\cite[\S~6]{Cu-Li}): Let $K$ be a number field with roots of unity
$\mu = \gekl{\pm 1}$. Choose a free abelian subgroup $\Gamma$ of $K\reg$
such that $K\reg = \mu \times \Gamma$. We obtain for the K-theory of the
ring C*-algebra $\fA[\so]$ attached to the ring of integers $\so$ of $K$:
\bgloz
  K_*(\fA[\so]) \cong
  \bfa
    K_0(C^*(\mu)) \otimes_{\Zz} \extalg (\Gamma) & \falls \# \gekl{v_{\Rz}} = 0, \\
    \extalg (\Gamma) & \falls \# \gekl{v_{\Rz}} \geq 1.
  \efa
\egloz
The distinction between the formulas in the two different cases
corresponds to a natural identification on the level of generators. As
abstract groups one obtains the same K-theory independently of the
number of real embeddings.

\section{The correct computation}

Let us first of all explain what went wrong in our original computation in \cite[\S~6.4, case~c)]{Cu-Li}: Let $\theta \in \Aut(C_0(\Rz))$ be the flip, i.e. $\theta (f) (x) = f(-x)$ for all $f \in C_0(\Rz)$. By equivariant Bott periodicity, we know that
$$
K_i(C_0(\Rz^2) \rtimes_{\theta \otimes \theta} \Zz / 2 \Zz) \cong
\bfa
  \Zz^2 & \falls i=0, \\
  \gekl{0} & \falls i=1.
\efa
$$
In the first part of the proof of \cite[Lemma~6.4]{Cu-Li}, we have claimed
that the automorphism $\id \otimes \theta$ of $C_0(\Rz^2) \rtimes_{\theta
\otimes \theta} \Zz / 2 \Zz$ acts as
$
\rukl{
\begin{smallmatrix}
  1&0 \\
  0&-1
\end{smallmatrix}
} $ in K-theory (in \cite{Cu-Li}, $\id \otimes \theta$ is denoted by
$\hat{\beta}_{(1,-1)}$). This however cannot be true. The reason is
that using the Pimsner-Voiculescu sequence, we would obtain as an
immediate consequence that $K_0(C_0(\Rz^2) \rtimes_{\theta \otimes
\theta} \Zz / 2 \Zz \rtimes_{\id \otimes \theta} \Zz) \cong \Zz
\oplus (\Zz / 2 \Zz)$. But as Lemma~\ref{correct-lem} below shows,
the correct result is $K_0(C_0(\Rz^2) \rtimes_{\theta \otimes
\theta} \Zz / 2 \Zz \rtimes_{\id \otimes \theta} \Zz) \cong \Zz$.

In the first part of the proof of \cite[Lemma~6.4]{Cu-Li}, we had considered
the number field $K = \Qz[\sqrt{2}]$ with ring of integers $\so = \Zz + \Zz \sqrt{2}$. The problem in our original computation was that we have assumed that in this particular case, the element $[u^1]_1 \times [u^{\sqrt{2}}]_1$ is part of a $\Zz$-basis for $G_{inf} \subseteq K_0(C^*(\so \rtimes \mu))$ (in the terminology of \cite[Lemma~6.1]{Cu-Li}). But this is not the case, only up to finite index. This is why \cite[Lemma~6.4]{Cu-Li} is false.

Here is now the correct computation:
\blemma
\label{correct-lem}
$
K_i(C_0(\Rz^2) \rtimes_{\theta \otimes \theta} \Zz / 2 \Zz \rtimes_{\id \otimes \theta} \Zz) \cong \Zz
$
for $i=0,1$.
\elemma
\bproof
The first step is the following simple observation:
\bgln
\label{isom}
  && C_0(\Rz^2) \rtimes_{\theta \otimes \theta} \Zz / 2 \Zz \rtimes_{\id \otimes \theta} \Zz / 2 \Zz \\
  &\cong& (C_0(\Rz) \otimes C_0(\Rz)) \rtimes_{\theta \otimes \theta} \Zz / 2 \Zz \rtimes_{\id \otimes \theta} \Zz / 2 \Zz
  \nonumber \\
  &\cong& (C_0(\Rz) \otimes C_0(\Rz)) \rtimes_{\theta \otimes \id} \Zz / 2 \Zz \rtimes_{\id \otimes \theta} \Zz / 2 \Zz \nonumber \\
  &\cong& ((C_0(\Rz) \rtimes_{\theta} \Zz / 2 \Zz) \otimes C_0(\Rz)) \rtimes_{\id \otimes \theta} \Zz / 2 \Zz \nonumber \\
  &\cong& \eckl{C_0(\Rz) \rtimes_{\theta} \Zz / 2 \Zz} \otimes \eckl{C_0(\Rz) \rtimes_{\theta} \Zz / 2 \Zz}. \nonumber
\egln
To get from the second to the third line, we just made use of the automorphism $(\Zz / 2 \Zz)^2 \cong (\Zz / 2 \Zz)^2$ given by $t_1 \ma t_1 t_2$, $t_2 \ma t_2$. Here $t_1$ and $t_2$ are the generators of the two copies of $\Zz / 2 \Zz$.

Since $K_0(C_0(\Rz) \rtimes_{\theta} \Zz / 2 \Zz) \cong \Zz$ and $K_1(C_0(\Rz) \rtimes_{\theta} \Zz / 2 \Zz) \cong \gekl{0}$ (see \cite[\S~3.3, Equation~(12)]{Cu-Li}), we deduce
\bgl
\label{K}
  K_i(C_0(\Rz^2) \rtimes_{\theta \otimes \theta} \Zz / 2 \Zz \rtimes_{\id \otimes \theta} \Zz / 2 \Zz) \cong
  \bfa
  \Zz & \falls i=0, \\
  \gekl{0} & \falls i=1.
  \efa
\egl

Now consider the automorphism $(\id \otimes \theta)\hat{\text{ }}$ of $C_0(\Rz^2) \rtimes_{\theta \otimes \theta} \Zz / 2 \Zz \rtimes_{\id \otimes \theta} \Zz / 2 \Zz$ which is dual to the action of the second copy of $\Zz / 2 \Zz$. Under the isomorphism~\eqref{isom}, $(\id \otimes \theta)\hat{\text{ }}$ corresponds to the automorphism $\hat{\theta} \otimes \hat{\theta}$, where $\hat{\theta}$ is the automorphism on $C_0(\Rz) \rtimes_{\theta} \Zz / 2 \Zz$ dual to $\theta$. Since $\hat{\theta}$ is either $\id$ or $-\id$ on $K_0(C_0(\Rz) \rtimes_{\theta} \Zz / 2 \Zz) \cong \Zz$, we conclude that
\bgl
\label{dual}
  ((\id \otimes \theta)\hat{\text{ }})_* = \id \text{ on }
  K_0(C_0(\Rz^2) \rtimes_{\theta \otimes \theta} \Zz / 2 \Zz \rtimes_{\id \otimes \theta} \Zz / 2 \Zz) \cong \Zz.
\egl

Plugging \eqref{K} and \eqref{dual} into the exact sequence from
\cite[Theorem~10.7.1]{Bla}, which connects the $K$-theory of the
crossed products by $\Zz$ and by $\Zz/2$ induced by $\id \otimes
\theta$ respectively, we obtain \bgloz
  K_i(C_0(\Rz^2) \rtimes_{\theta \otimes \theta} \Zz / 2 \Zz \rtimes_{\id \otimes \theta} \Zz) \cong \Zz \text{ for } i=0,1.
\egloz
\eproof

With this lemma, the computation of the K-theory of the ring C*-algebras follows the same line of arguments as in \cite{Cu-Li}. Let us explain this briefly using the same notations as in the introduction and as in \cite[\S~6.4, case~c)]{Cu-Li}. Combining \cite[(4)]{Cu-Li} with \cite[Corollary~4.2]{Cu-Li} and using a refined version of \cite[Lemma~6.3]{Cu-Li}, it is straightforward to see that the K-theory of $\fA[\so]$ coincides with the K-theory of $C_0(\Az_{\infty}) \rtimes K\reg$. As in \cite[\S~6.4, case~c)]{Cu-Li}, let $K \reg = \mu \times \Gamma$ and choose a $\Zz$-basis $\gekl{p, p_1, p_2, \dotsc}$ of $\Gamma$, with $p \in \Zz \pos$. We can arrange that $\# \menge{v_{\Rz}}{v_{\Rz}(p_1)<0}$ is odd and $\# \menge{v_{\Rz}}{v_{\Rz}(p_i)<0}$ is even for all $i>1$. Let $\Gamma_m = \spkl{p, \dotsc, p_m}$ and $\Gamma'_m = \spkl{p, p_2 \dotsc, p_m}$. An iterative application of the Pimsner-Voiculescu sequence gives
\bgloz
  K_*(C_0(\Az_{\infty}) \rtimes (\mu \times \Gamma_m)) \cong \extalg (\Gamma'_m)
\egloz
and thus
\bgloz
  K_*(\fA[\so]) \cong \extalg (\Gamma).
\egloz

\end{document}